\documentclass[extra, camera, onecolumn,   
]{gji}
\usepackage{timet}

\usepackage[utf8]{inputenc}
\usepackage{lmodern}
\usepackage[T1]{fontenc}
\usepackage{amssymb, amsmath, amsfonts, multirow, subcaption, booktabs}
\usepackage{booktabs}
\usepackage{xspace}
\usepackage[final]{graphicx}
\usepackage{ifthen}
\usepackage[addedmarkup=bf]{changes}
\usepackage[final, allcolors=blue, colorlinks=true]{hyperref}
\usepackage[modulo]{lineno}

\makeatletter
\let\zz@tabular\@tabular
\let\zzendtabular\endtabular
\let\zz@xtabularcr\@xtabularcr
\let\zz@tabclassz\@tabclassz
\let\zz@tabclassiv \@tabclassiv 
\let\zz@tabarray\@tabarray
\makeatother

\newlength{\fwidth}
\newlength{\cwidth}
\newlength{\twidth}
\ifthenelse{\boolean{GJ@referee}}{%
  \setlength{\cwidth}{234pt}%
  \setlength{\fwidth}{468pt}%
  \setlength{\twidth}{468pt}}{%
  \setlength{\fwidth}{504pt}%
  \ifthenelse{\boolean{@twocolumn}}{%
    \setlength{\cwidth}{240pt}%
    \setlength{\twidth}{240pt}}{%
    \setlength{\cwidth}{252pt}%
    \setlength{\twidth}{504pt}%
  }
}
\graphicspath{{./figures/}}

\definecolor{alizarin}{rgb}{0.82, 0.1, 0.26}
\definecolor{mygreen}{rgb}{0, .5, 0}
\setaddedmarkup{\color{mygreen}{\uline{#1}}}
\setdeletedmarkup{\color{alizarin}{\uwave{#1}}}
\definechangesauthor[color=blue]{aamor}

\usepackage{siunitx}
\usepackage{tabularx}

\AtBeginDocument{
  \label{CorrectFirstPageLabel}
  
}

\newcommand{\petgem}{\texttt{PETGEM}\xspace}
\newcommand{\empymod}{\texttt{empymod}\xspace}

\title[Tailored meshing for parallel 3D CSEM using HEFEM]{Tailored meshing for parallel 3D electromagnetic modeling using high-order edge elements\footnote{Journal: Journal of Computational Science\\ Corresponding author: Octavio Castillo-Reyes\\ First author: Octavio Castillo-Reyes\\ Received at Editorial Office: 15 Dec 2021}}
\author[O. Castillo-Reyes \emph{et al.}]
  {\Large Octavio Castillo-Reyes$^1$,   
   Adrian Amor-Martin$^2$,                    
   Arnaud Botella$^3$,                    
   Pierre Anquez$^3$, \\                    
   \Large Luis Emilio Garc\'ia-Castillo$^2$\\                    
    \footnotesize
  $^1$ Barcelona Supercomputing Center (BSC), Plaça Eusebi Güell 1-3, 08034 Barcelona (ES). E-mail: \href{mailto:octavio.castillo@bsc.es}{octavio.castillo@bsc.es}\\[-.3em]
   \footnotesize
  $^2$ Department of Signal Theory and Communications. University Carlos III of Madrid, 28903 Madrid (ES) \\[-.3em]
  \footnotesize
  $^3$ Geode-solutions, Technopole H\'elioparc 2 Av. du Pr\'esident Pierre Angot, 64000 Pau (FR) \\[-.3em]
}

\pagerange{\pageref{firstpage}--\pageref{lastpage}}
\volume{}
\pubyear{2021}

\begin{document}

\label{firstpage}

{\makeatletter
\let\@tabular\zz@tabular
\let\endtabular\zzendtabular
\let\@xtabularcr\zz@xtabularcr
\let\@tabclassz\zz@tabclassz
\let\@tabclassiv \zz@tabclassiv 
\let\@tabarray\zz@tabarray
\maketitle
}

\begin{summary}
We present numerical experiments for geophysics electromagnetic (EM) modeling based upon high-order edge elements and supervised $h+p$ refinement approaches on massively parallel computers. Our high-order $h+p$ refinement strategy is based on and extends the \petgem code. We focus on the performance study in terms of accuracy, convergence rate, and computational effort to solve real-life 3D setups based on synthetic and experimental data for energy reservoir characterization. These test cases show variable resolution discretization needs and realistic physical parameters. In general, our numerical results are consistent theoretically. The use of $h-$adapted meshes was efficient to achieve a certain accuracy level in the synthetic EM responses. Regarding global $p-$refinement, $p=2$ exhibits the best accuracy/performance trade-off. Selective $p$-refinement might offer a better compromise between accuracy and computational cost. However, for $p-$refinement at different entities, the best refinement scheme consists of using $p=3$ at the volume level with $p=1$ at faces and edges. Thus, $p-$refinement can be competitive if applied hierarchically. Nevertheless, we acknowledge that the performance of our supervised $h+p$ refinement strategy depends on the input model (e.g., conductivity, frequency, domain decomposition strategy, among others). Whatever the chosen configuration, our numerical results provide an in-depth understanding of EM modeling's pros and cons when supervised $h+p$ refinement schemes are applied.
\end{summary}


\section{Introduction}
\label{introduction}
Electromagnetic (EM) modeling routines play a key role in studying and interpreting subsurface electric conductivity distribution. They are widely used in academia and industry because of their capacity to reduce ambiguities in interpreting geophysical datasets through mapping conductivity variations in the subsurface. As a result, there are abundant  schemes to implement 3D EM forward modeling algorithms (see \citet{SG.05.Avdeev,SG.10.Borner} for a detailed review). For inversion of EM datasets, a large number of forward modeling computations are needed. Therefore, 3D EM modeling (and inversion) algorithms should be particularly sought for:
\begin{enumerate}
    \item Providing accurate solutions in a feasible run-time, although the uncertainties associated with the domain discretization, numerical operator, among others.
    \item Tackling problems efficiently using cutting-edge computing architectures, including high-performance computing (HPC), thus ability to deal with real-life models.
    \item Modularity and flexibility to cope with a variety of real-life model workflows with the possibility to easily add or remove components without having to rewrite large parts of the code.
\end{enumerate}
For 3D EM modeling, multiple references use N\'ed\'elec FE of low-order polynomial functions (e.g., first-order)~\citep{SG.10.Borner,Farquharson2011,Mukherjee2011,Ren2013,Um2013,Cai2014,Chung2014,Grayver2014,Kordy2015,Castillo2017,CAG.18.CastilloReyes}, and high-order polynomial functions~\citep{Schwarzbach2011,Grayver2015,GJI.19.CastilloReyes,Rochlitz2019}. N\'ed\'elec FE offer a proper mechanism to discretize the space, $\mathbf{H}$(\textit{curl}), of the electric vector field. 
Also, it is possible to apply adaptive mesh refinement to improve the solution accuracy with decreased computational effort compared to uniform refinement. The mesh refinement is an elaborated method in FE modeling aiming at designing locally refined meshes to accurately model the fields in specific regions of interest while using as few mesh elements as possible. Two common refinement methods of achieving more accurate FE solutions are to increase the number of elements ($h-$refinement) and to employ higher-degree interpolation functions ($p-$refinement). Furthermore, sophisticated methods for a posteriori error estimation combined with both refinement strategies also referred to as $hp-$refinement, might offer exponential convergence rates~\citep{Grayver2015}. These mesh refinement strategies can be automatic or non-automatic.

Adaptive mesh $h-$refinement has already been investigated in geophysical electromagnetics~\citep{Plessix2007,Schwarzbach2011,CAG.18.CastilloReyes}. However, there is a lack of studies about high-order N\'ed\'elec finite elements (HEFEM) and their impact on 3D EM modeling in the realm of 3D geo-electromagnetic modeling (models that include experimental data and not only synthetic). The exceptions are the works performed by~\cite{Schwarzbach2011} and~\cite{Grayver2015}, which demonstrate the advantages of high-order basis functions regarding the needs of degrees of freedom (dof) to satisfy the prior chosen quality criteria in the EM responses. HEFEM have also been shown to be beneficial in other areas of EM wave modeling, such as cavity analysis~\citep{jian2003,luise_cmame11_hofeim}, brain microwave studies~\citep{bonazzoli_ijnm18}, among others with smooth solutions~\citep{bergot_jcp13,olm_aes19,Eisentrager2020}.

Our contribution is twofold. First, we introduce a robust numerical scheme that can be applied to simulate different and relevant scenarios in the field of the 3D Controlled-source Electromagnetic Method (CSEM). Second, we provide a firm basis to justify the development of more robust techniques such as a fully automatic goal-oriented approach~\citep{Key2011,Pardo2011,Schwarzbach2011,Grayver2015}.
Hence, the overall research hopes to contribute by providing a comprehensive and quantitative study about the use of HEFEM in conjunction with unstructured and $h+p$ adapted meshes for real-scale geo-electromagnetic modeling. We acknowledge that high-order polynomial interpolation functions for EM modeling in geophysics are previously reported by \citep{Schwarzbach2011,Grayver2015,GJI.19.CastilloReyes,Rochlitz2019}. 
The authors of these studies investigate convergence and numerical accuracy for polynomial interpolation functions of order $p=1,2,3$. However, our research provides relevant information for an in-depth understanding of the pros and cons when still higher polynomial schemes ($p=1,2,3,4,5,6$) and supervised $h+p$ adapted meshing are employed. We state that our experiments and conclusions can be considered a preliminary stage that favors simplicity and practicality, which can be helpful not only for experts but also for the general community interested in EM modeling. In this way, we set a bound to apply automatic refinements that will be explored in future works. Finally, to verify the robustness of our numerical strategies, we simulate 3D CSEM setups with relevance for both academia and industry. Each 3D CSEM setup under consideration presents a particular numerical modeling challenge, being a suitable approach to studying the proposed numerical schemes. Furthermore, we stress that we not only analyze models with synthetic data but also models with experimental data, a clear contribution with respect to the rest of the state-of-the-art works.

On top of that, the core motivations for this study are to evaluate the benefits and limitations of HEFEM and supervised $h+p$ tailored meshing for real-life 3D CSEM surveys. We focus on investigating its performance in terms of convergence rate, CPU time, and memory requirements. We use a supervised approach, which consists of:

\begin{enumerate}
    \item Determining prior rules to build $h-$adapted meshes for a set of basis orders $p=1,2,3,4,5,6$.
    \item Applying a high-order basis discretization on regions (known a priori) where the solution is not smooth (e.g., transmitter vicinity). To increase the $p$ order on specific regions without sacrificing computational performance, we consider a constant number of dof per element.
\end{enumerate}
To compute synthetic EM responses, we have used the \petgem code~\citep{CAG.18.CastilloReyes}, which has been shown as an efficient modeling routine on massively parallel computers.  We point out that our numerical scheme and its results are independent of the electric field formulation (e.g. total or secondary field formulation), transmitter type (e.g. electric or magnetic dipole), transmitter signal (e.g. time-harmonic, direct current or transient), and conductivity material properties (e.g. isotropic or anisotropic). Still, to preserve brevity, we particularized our study to 3D CSEM modeling for marine and land applications in the context of energy reservoirs characterization (oil $\&$ gas, and geothermal) based on synthetic and experimental data.

The rest of the paper is organized as follows. Section~\ref{theory} introduces the governing equations for the EM modeling under consideration and its discretization using the HEFEM. Also, we provide details about our supervised $h+p$ refinement scheme. In Section~\ref{numerical_validation}, we perform numerical simulations to investigate the numerical scheme's robustness for different high-order FE basis. Also, we extensively discuss essential points that control the computational performance and suitability of our supervised adaptive mesh refinement technique in geophysical electromagnetics. Finally, Section~\ref{conclusions} provides summary remarks and conclusions.
\section{Theory}
\label{theory}

\subsection{Governing equations}
\label{governing_equations}
The EM phenomena under consideration is described mathematically by the frequency-domain Maxwell's equations in a diffusive form. Based on this approach, displacement currents are neglected. By assuming a time-harmonic dependence expressed by $e^{-i \omega t}$, these equations can be written as
\begin{align}
	\boldsymbol{\nabla} \times \mathbf{E} &= i\omega \mu_{0}\mathbf{H}, \label{eq:maxwell_diffusive_form1} \\
	\boldsymbol{\nabla} \times \mathbf{H} &= \mathbf{J}_{\text{s}} + \sigma \mathbf{E}, \label{eq:maxwell_diffusive_form2}
\end{align}
where $\mathbf{E}$ is the electric field, $\mathbf{H}$ is the magnetic field, $i$ is the imaginary unit, $\omega$ is the angular frequency, $\mu_{0}$ is the free-space magnetic permeability, $\mathbf{J}_{\text{s}}$ is the distribution of source current, $\sigma \mathbf{E}$ is the induced current in the conductive Earth, and $\sigma$ is the electric conductivity tensor.

It is frequently convenient to have a formulation in terms of the total field. After substituting eq.~\eqref{eq:maxwell_diffusive_form1} into eq.~\eqref{eq:maxwell_diffusive_form2}, we obtain
\begin{align}
\boldsymbol{\nabla} \times \boldsymbol{\nabla} \times \mathbf{E} - i \omega \mu_{0} \sigma \mathbf E = i \omega \mu_{0} \mathbf{J}_{\text{s}},
\label{eq:electric_field_weak}
\end{align}
which is also known as the \textit{curl-curl} formulation of the problem in terms of the total field. This approach can avoid numerical errors that arise when the source is located within the region of anomalous properties (e.g. models with high conductivity contrasts or with bathymetry/topography variations). The disadvantage is that a slightly larger computational domain is required in order to discretize the source properly and avoid artifacts arising from the reflections on the artificial boundary conditions of the domain.

\subsection{HEFEM Solution}
\label{hefem}
For the numerical solution of eq.~\eqref{eq:electric_field_weak}, we consider a 3D computational domain discretized in a set of tetrahedral finite elements. For a comprehensive introduction to the FE method, we refer to~\cite{luise_book98,jian2003,monk_book03}, and~\cite{Thompson2005}. This section only provides a brief outline for understanding the HEFEM implemented within the \petgem code.

The space that we have chosen for the electric field $\mathbf E$ in
the domain $\Omega$ is the \textit{curl}-conforming space 
\begin{align}
\mathbf{H}(\text{curl},\Omega) := \{ \mathbf{w} \in \left[ L_2(\Omega) \right]^3 \; \mid \; \nabla\times\mathbf{w} \in \left[ L_2(\Omega) \right]^3 \}
\label{eq:hcurl}
\end{align}
with $L_2(\Omega)$ as the space of square integrable functions in $\Omega$. Specifically, we use
\begin{equation}
\begin{aligned}
    \mathbf{H}_0(\text{curl},\Omega) := \{ & \mathbf{w} \in \mathbf{H}(\text{curl},\Omega) \mid \mathbf{n} \times \mathbf{w} = \mathbf{0}\, \text{on} \,\partial\Omega_{\text{D}} \},
    \label{eq:hzcurl}
\end{aligned}
\end{equation}
within the domain $\Omega$ subject to a homogeneous Dirichlet boundary condition $\mathbf{n} \times \mathbf E = \mathbf 0$ on the domain boundary $\partial\Omega$. 

To discretize the eq.~\eqref{eq:electric_field_weak},
we use the hierarchical basis functions by~\cite{fuentes_camwa15}, that belong to the mixed-order family proposed by~\cite{nedelec80} and span
a \textit{curl}-conforming space $\mathcal{W}_k$. To allow non-uniform $p-$refinement, we use the hierarchical property of the basis functions (e.g., we use incremental spaces of 
basis functions $\tilde{\mathcal{W}_k}$), so the space $\mathcal{W}_k$
is
\begin{align}
    \mathcal{W}_k = \tilde{\mathcal{W}_1} \oplus \tilde{\mathcal{W}_2} \oplus \dots \oplus \tilde{\mathcal{W}_k}.
    \label{eq:basis_functions}
\end{align}
The application of the Galerkin procedure to eq.~\eqref{eq:electric_field_weak} leads to the variational formulation $\mathbf{E} \in \mathbf{H}_0(\text{curl},\Omega)$ such that
\begin{align}
    \iiint \limits_{\Omega} (\nabla\times\mathbf{W}) \cdot (\nabla\times\mathbf{E}) d\Omega - i\omega\mu_{0} \iiint \limits_{\Omega} \mathbf{W} \cdot \sigma\mathbf{E} d\Omega &=  
    i\omega\mu_{0}\iiint \limits_{\Omega} \mathbf{W} \cdot \mathbf{J}_{\text{s}} d\Omega \quad \forall \,\mathbf{W} \in \mathbf{H}_{0}(\text{curl},\Omega).
  \label{eq:var_form}
\end{align}
The discretization of this variational problem provides a 
system of equations $\mathbf{A} \mathbf{x} = \mathbf{b}$, 
where the vector of unknowns $\mathbf{x}$ allows us to approximate 
the electric field in any point of the computational domain by using 
\begin{align}
    \tilde{\mathbf{E}} = \sum_{i=1}^{N_e} x_i \mathbf{w}_i,
    \label{eq:interpolation_basis}
\end{align}
where $\tilde{\mathbf{E}}$ is the approximation of the electric field, 
$N_e$ the size of the matrix to be solved, $x_i$ are the weights of the basis functions obtained from $\mathbf{A} \mathbf{x} = \mathbf{b}$, and $\mathbf{w}_i$ the basis functions used to solve the FE problem.

\subsection{Supervised \textit{h+p} adaptive meshing}
\label{hp_refinement}
In our previous study~\cite{GJI.19.CastilloReyes} we have introduced a methodology to design $h-$adapted tetrahedral meshes for the CSEM problem. We used the skin-depth ($\delta$) principle as the main quality criteria to determine characteristic mesh sizes for piecewise orders $p=1,2,3$. In this paper, we focus on extending this strategy to higher-order variants ($p=4,5,6$) on computational domains with general anisotropy. Thus, we provide a comprehensive description for the piecewise family implemented within our modeling routine. The procedure, which is also guided by both the physical parameters and piecewise order, is composed of two steps:
\begin{enumerate}
    \item Computation of the characteristic spacing $d_{\delta}$, which can be expressed as
    \begin{align}
    d_{\delta}(f, p) = \frac{\delta_{\min}(f)}{\lambda_{\delta}(p)},
    \label{eq:skin_depth_frequency}
    \end{align}
    where $f$ is the modeling frequency, $p$ is the piecewise order, $\delta_{\min}$ is the minimum skin depth in the input conductivity/resistiviy model, and $\lambda_{\delta}$ is the number of points per skin depth. 
    \item Computation of local refinement $d_{s}$ that is applied close to regions of interest such as source and receiver locations. The formal expression of $d_{s}$ can be written as
    \begin{align}
    d_{\text{s}} = \min\left(\frac{L_{\text{s}}}{r_{\text{s}}(p)}, d_{\delta} (f, p)\right),
    \label{eq:source_refinement}
    \end{align}
    where $L_{\text{s}}$ is the source dipole length and $r_{\text{s}}$ is a resolution number that also depends on the piecewise order $p$. 
\end{enumerate}
Based on this strategy, all cell sizes in the computational domain are constrained by rules that satisfy the prior chosen quality criteria in the EM responses. In Section~\ref{test_layered_meshing_rules} of this paper, we follow a rigorous methodology to determine $\lambda_{\delta}$ and $r_{\text{s}}$. 
\begin{table}
\caption{Number of dof per tetrahedral element.}
\begin{center}
\begin{tabular}{c c c}
\toprule
\multicolumn{2}{c}{Tetrahedron}             & $\frac{1}{2}p(p+2)(p+3)$  \\
\cmidrule{1-3}
    & Edge                                  & $p$                      \\     
    & Face                                  & $p(p-1)$                 \\   
    & Volume                                & $\frac{1}{2}p(p-1)(p-2)$ \\
\bottomrule   
\end{tabular}
\end{center}
\label{table:dof}
\end{table}

The implemented hierarchical basis functions are based on orthogonal polynomials (Jacobi for face and interior functions, Legendre for edge) which allow that each entity (edge, face or interior functions) have a different approximation order within the same element. For convenience, Table~\ref{table:dof} shows the number of dof per tetrahedral element within the computational mesh.

\subsection{Parallel implementation}
The supervised adaptive mesh refinement scheme described in the previous section has been implemented in a fully distributed fashion (MPI standard) using Python language and \texttt{mpi4py}\xspace~\citep{Dalcin2008}. Our implementation is based on and extends the \petgem code. The algorithm supports general anisotropic cases and takes advantage of the inherent features within \petgem for parallel computations. This section only provides a brief outline for understanding the HPC implementation. For a comprehensive introduction to the \petgem code and its computational performance, we refer to~\citep{CAG.18.CastilloReyes,GJI.19.CastilloReyes,TGRS.21.CastilloReyes,CAG.22.CastilloReyes}. 

The \petgem workflow is composed of four main phases: \texttt{pre-processing}, \texttt{assembly}, \texttt{solver}, and \texttt{post-processing}. The first one is responsible for providing functions to change input data related to mesh ($h-$refinement) into a more suitable representation for \petgem. These are, in particular, the nodal spatial coordinates, the mesh connectivity, and boundary conditions. Further, additional information needed for handling matrix storage allocation is retrieved here. The \texttt{assembly} phase encapsulates information and functions coming from the HEFEM and $p-$refinement. This phase particularly comprises the computation of the \textit{curl}-conforming basis functions and the computation of the elemental FE matrices. Further, in the \texttt{solver} phase, the system of linear equations is solved directly using \texttt{MUMPS}\xspace~\citep{MUMPS2006} or iteratively via the \texttt{petsc4py}\xspace package~\citep{Dalcin2011}. These libraries offer a large selection of parallel iterative Krylov solvers and multifrontal direct methods. Finally, in the \texttt{post-processing} stage, the EM responses are computed.

\begin{figure}
	\centering
	\includegraphics[trim= 0 0 0 0,clip,width=0.7\textwidth]{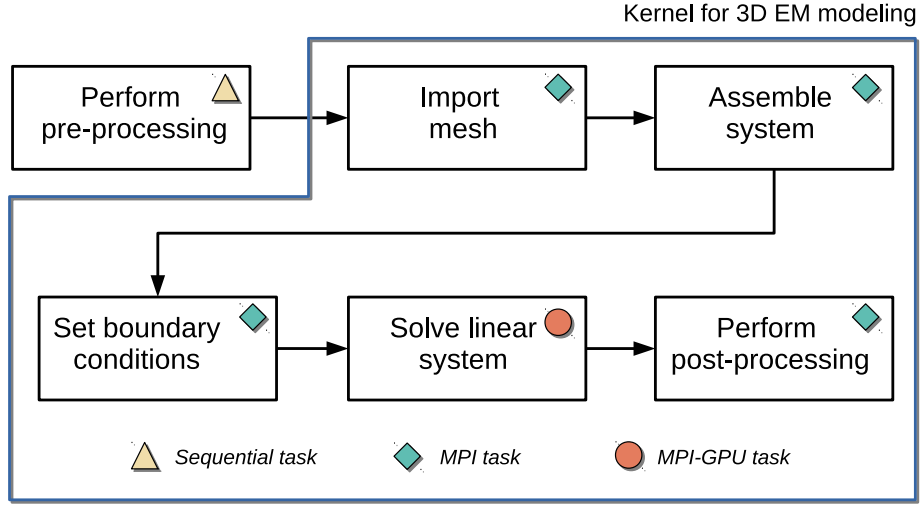}
	\caption{Outline of the overall \petgem work-flow. The \texttt{pre-processing} phase is sequential while the remaining stages are parallel. Notice that the \texttt{solver} phase has support for both MPI and MPI-GPUs platforms. Such an HPC implementation can be useful and transparent in the absence of GPUs devices.}
	\label{fig:petgem_work_flow}
\end{figure}

The stages mentioned above are used to fuels a \petgem \texttt{kernel} responsible for solving the input 3D CSEM model. The corresponding work-flow is depicted in Figure~\ref{fig:petgem_work_flow}.

\section{Numerical Validation}
\label{numerical_validation}
This section describes a set of numerical experiments with HEFEM and supervised $h+p$ refinement to solve realistic and challenging 3D CSEM setups. The numerical results presented in this section are divided into two sections. First, we perform a series of tests to verify the implementation of HEFEM. This set of tests corresponds to the so-called Method of Manufactured Solutions (MMS)~\citep{marchand_ap14,garcia_apm16}. Second, we perform 3D CSEM simulations on different setups to find out which basis order offers the best compromise between accuracy and computational effort.

The process of verification and analysis of our meshing strategy consists of three steps:
\begin{enumerate}
    \item Robustness study of uniform $p-$refinement (experiment~\ref{test_mms})
    \item Determine the optimal element size (uniform $h-$refinement) for each basis order (experiment~\ref{test_layered_meshing_rules})
    \item Determine the appropriate combination of basis order at different element levels (experiment~\ref{marine_mr3d_model})
    \item Particularize the best configuration obtained from the previous points to regions where high numerical accuracy is required (experiment~\ref{land_valles_model})
\end{enumerate}
We point out that although our modeling tool supports 3D anisotropy resistivity, we only consider models with Vertical Transverse Isotropy (VTI). Furthermore, we state 
that our meshing rules do not consider anisotropy (an in-depth analysis of this effect is beyond the scope of this paper). However, this does not have a negative impact on our numerical experiments for two reasons. First, because of the diffusive nature of the problem (the EM field decays exponentially). Second, the modeling configurations require that the EM field be measured at sites distant from the source location (several skin-depths away). To compute the $L^{2}$-norm errors and their corresponding convergence orders $\mathcal{O}_{L^{2}}$ of our numerical solutions, we use the expressions suggested by~\cite{GJI.19.CastilloReyes}. To perform the unstructured tetrahedral mesh generation, we employ the \texttt{Gmsh}\xspace tool~\citep{Geuzaine2008}. All simulations have been carried out on \textit{Marenostrum} supercomputer.

\subsection{MMS}
\label{test_mms}
We verify the $p-$refinement scheme on an MMS test as the first example. From eq.\eqref{eq:electric_field_weak}, we need to include non-homogeneous Dirichlet boundary conditions on the boundary of the problem. We manufacture an analytical solution to the differential equation (by solving the problem backward), and then we measure the error between the approximated electric field and the manufactured solution~\citep{marchand_ap14,garcia_apm16}.

The boundary value problem under consideration is expressed as
\begin{subequations}
    \begin{align}
        \nabla\times(\mu_r^{-1}\nabla\times \mathbf{E}) - i\omega\mu_0\sigma\mathbf{E} &= \mathbf{F} && \text{in}\ \Omega, \\
        \mathbf{n}\times(\mathbf{E}\times\mathbf{n}) &= \mathbf{\Psi}_{\text{D}} && \text{on}\ \Gamma_{\text{D}},
    \end{align}
    \label{eq:bvp}
\end{subequations}
with a source excitation defined as
\begin{subequations}
    \begin{align}
        \mathbf{F} &= \nabla\times\left(\nabla\times\mathbf{E}_\text{MMS}\right) - i\omega\mu_0\sigma\mathbf{E}_\text{MMS}, \\
        \mathbf{\Psi}_{\text{D}} &= \mathbf{n}\times(\mathbf{E}_\text{MMS}\times\mathbf{n}).
    \end{align}
\end{subequations}
To assign non-homogeneous Dirichlet boundary conditions to a hierarchical set of basis functions, we need to project the solution on the boundary faces~\citep{demkowicz_book06}. This field projection is performed in two steps. First, for each edge $e$ we solve 
\begin{align}
    \int \limits_{e} (\mathbf{W}_{e}\cdot\hat{\mathbf{\tau}}) (\mathbf{E} \cdot \hat{\mathbf{\tau}}) de = \int \limits_{e} (\mathbf{W}_{e}\cdot\hat{\mathbf{\tau}}) (\mathbf{E}_{\text{MMS}} \cdot \hat{\mathbf{\tau}}) de, \label{eq:dir_edge}
\end{align}
where $\hat{\mathbf{\tau}}$ is the edge unit vector. In eq.~\eqref{eq:dir_edge} we belonging to the restriction of $\mathbf{H}(\text{curl},\Omega)$ to the edge. Second, for each face $f$ we solve
\begin{equation}
\begin{aligned}
    \int \limits_{f} (\hat{\mathbf{n}}\times(\mathbf{W}_{f}\times\hat{\mathbf{n}}))\cdot(\hat{\mathbf{n}}\times(\mathbf{E}\times\hat{\mathbf{n}})) df = \int \limits_{f} (\hat{\mathbf{n}}\times(\mathbf{W}_{f}\times\hat{\mathbf{n}})) \cdot \mathbf{d}_f df,
    \label{eq:dir_face}
\end{aligned}
\end{equation}
where we belonging to the restriction of $\mathbf{H}(\text{curl},\Omega)$ to the face. In eq.~\eqref{eq:dir_face} the excitation $\mathbf{d}_f$ is obtained through
\begin{align}
    \mathbf{d}_f = \hat{\mathbf{n}}\times(\mathbf{E}_\text{MMS}\times\hat{\mathbf{n}} - \sum_{\substack{e\in f}} \mathbf{d}_e),
\end{align}
and the contribution for each edge $\mathbf{d}_e$ is 
\begin{align}
    \mathbf{d}_e = \sum_{i=1}^{n_e} x_i (\hat{\mathbf{n}}\times(\mathbf{W}_{e}\times\hat{\mathbf{n}})),
\end{align}
being $n_e$ the number of dof for each edge, and 
$x_i$ the $i$-th component of the solution obtained from eq.~\eqref{eq:dir_edge}.

The analytical solution is set to a plane wave as example of a
smooth function, which is expressed as
\begin{align}
    \mathbf{E}_\text{MMS} = \mathbf{E}_\text{pol}e^{-jk_0(\mathbf{k}_p\cdot\mathbf{r})}, \label{eq:e_exp_mms}
\end{align}
where $\mathbf{E}_\text{pol} = \hat{\mathbf{x}} + \hat{\mathbf{y}}$, and
$\mathbf{k}_p = \hat{\mathbf{z}}$, with $\hat{\mathbf{x}}$, $\hat{\mathbf{y}}$,
and $\hat{\mathbf{z}}$ as the unit vectors in the X, Y, and Z axes.
This function is not a solution for eq.~\eqref{eq:var_form}, but it may serve for verification purposes. Our computational domain is a $[0,1]^3$ m cube in vacuum space. 
\begin{table}
\caption{Mesh hierarchies for a MMS test using uniform and non-uniform refinement in $p$. Hierarchical mesh level, number of tetrahedral elements and number of dof are provided.}
\begin{center}
\begin{tabular}{c c c c c c c c }
\toprule
\multirow{3}{*}{Mesh level} & \multirow{3}{*}{Elements} & \multicolumn{3}{c}{uniform refinement} & \multicolumn{3}{c}{non-uniform refinement} \\
\cmidrule{3-8}
&  & $p=1$ & $p=2$ & $p=3$ & $p\in[1,2]$ & $p\in[2,3]$ & $p\in[1,3]$ \\
\cmidrule{1-8}
1 & $727$   & $1\,163$    & $5\,636$   & $15\,600$ & $1\,163$ & $5\,636$ & $1\,163$  \\
2 & $4\,716$  & $6\,624$  & $33\,590$  & $95\,046$ & $16\,245$ & $53\,418$ & $36\,073$  \\
3 & $16\,028$ & $21\,093$ & $109\,562$ & $313\,491$ & $78\,198$  & $232\,091$  & $200\,727$\\
4 & $36\,395$ & $46\,546$ & $244\,340$ & $702\,567$ & $215\,458$ & $621\,983$  & $593\,101$ \\
5 & $71\,551$ & $89\,716$ & $474\,488$ & $1\,368\,969$ & $474\,488$ & $1\,368\,969$ & $1\,368\,969$\\
\bottomrule   
\end{tabular}
\end{center}
\label{tab:mesh_hierarchies_MMS}
\end{table}
\subsubsection{Uniform refinement}
\label{uniform_refinement_mms}
\begin{figure}
  \centering
  \includegraphics[trim= 0 0 0 0, clip,width=0.55\textwidth]{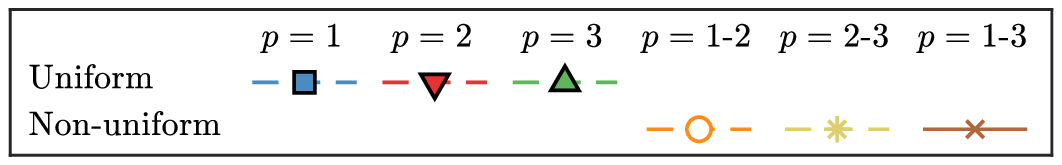}
  
  \begin{subfigure}[b]{0.45\textwidth}
    \includegraphics[trim= 0 0 0 0, clip,width=\textwidth]{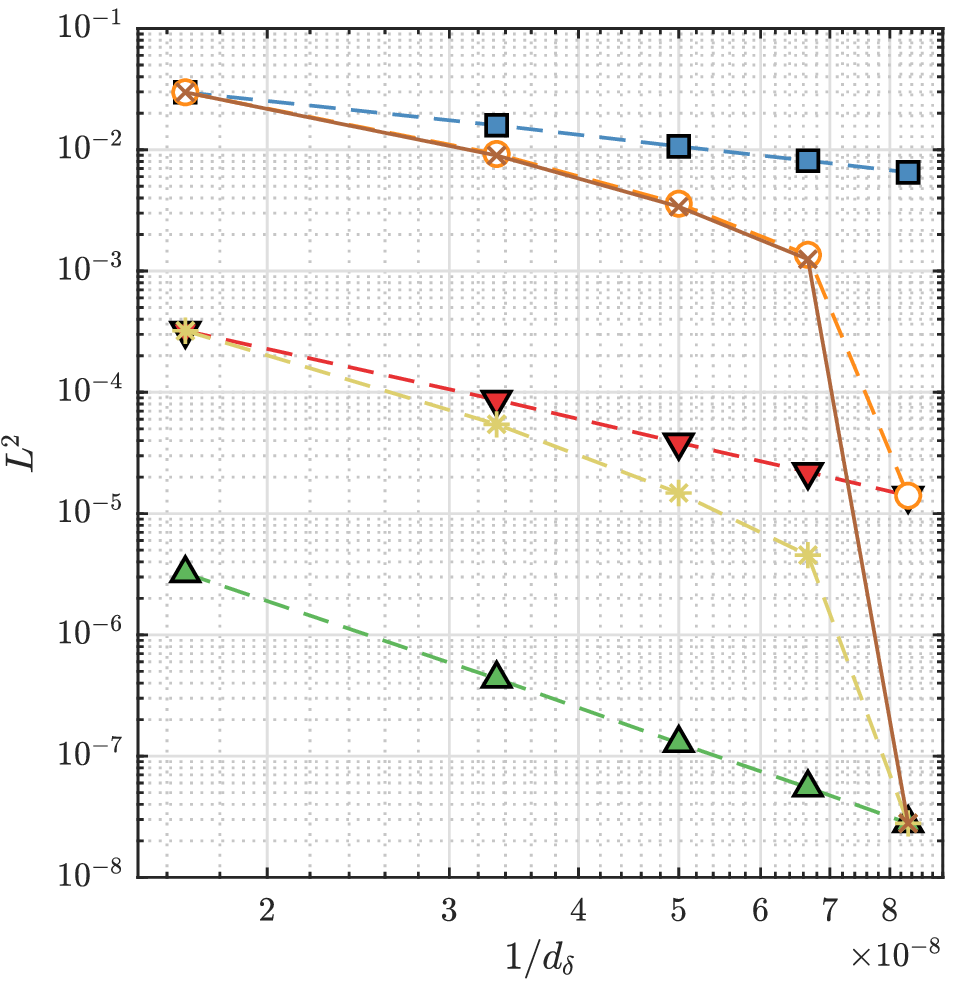}
    \caption{}
    \label{fig:convergence_mms_a}
  \end{subfigure}
  \hspace{.2cm}
  \begin{subfigure}[b]{0.45\textwidth}
    \includegraphics[trim= 0 0 0 0, clip, width=\textwidth]{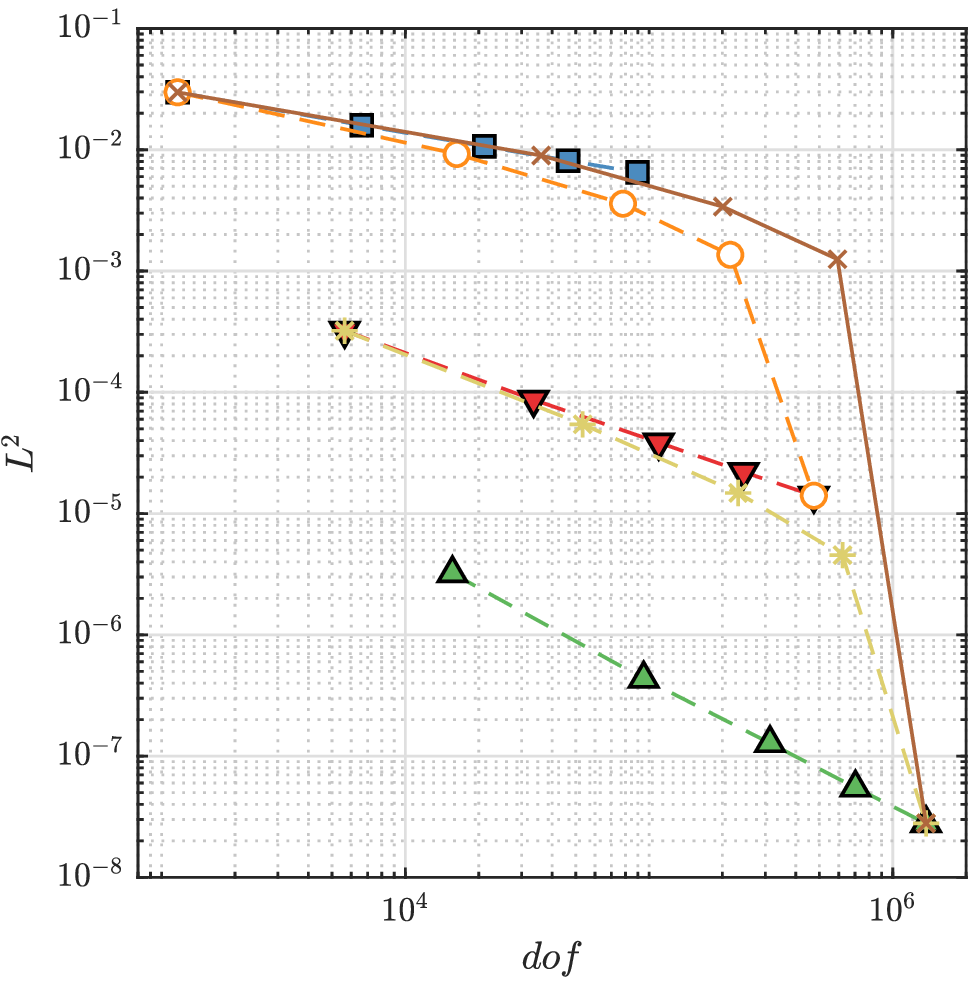}
    \caption{}
    \label{fig:convergence_mms_b}
  \end{subfigure}
  \caption{Convergence behavior of the plane wave with MMS in terms of (a) characteristic mesh spacing $d_\delta$ and (b) number of dof. Convergence results for uniform and non-uniform in $p$ are given.}
\end{figure}

First, numerical approximations have been computed on a set of globally $p-$refined meshes up to the piecewise $p=3$. The mesh hierarchies are given in Table~\ref{tab:mesh_hierarchies_MMS}. The convergence behavior regarding characteristic mesh size and number of dof is show in Fig.~\ref{fig:convergence_mms_a} and Fig.~\ref{fig:convergence_mms_b}, respectively. For piecewise $p=1,2,3$, excellent average slopes of $0.94$, $1.94$, and $2.95$ are observed. These results are consistent with the theoretical slopes $\mathcal{O}(h^p)$~\citep{luise_book98}. 
\subsubsection{Non-uniform refinement}
\label{non_uniform_refinement_mms}
Next, we changed the polynomial order of a given percentage mesh elements to assess that the obtained error lies within the mesh's minimum and maximum polynomial orders. Thus, we start in the same first mesh level and, for each finer mesh, we increase the piecewise order (two cases: $p+1$ and $p+2$) in $25\%$ of the mesh elements in a cumulative way. The resulting mesh statistics of this approach are given in Table~\ref{tab:mesh_hierarchies_MMS}. 

From this experiment, we obtained three new convergence curves labeled as $p=1$-$2$, $2$-$3$, $1$-$3$ and also depicted in Fig.~\ref{fig:convergence_mms_a} and Fig.~\ref{fig:convergence_mms_b}. For all three cases, the start- and end-points of convergence curves are the same as that attained on the previous homogeneous $p-$refinement test. The obtained error is between the curves associated with the minimum and maximum polynomial order in the mesh, which can be considered a validation of our non-uniform $p-$refinement scheme. Besides, the difference of accuracy between the fourth level (75\% refined elements) against the fifth level (100\% refined elements) is considerable. This error difference is due to the nature of the plane wave used as the excitation source. More concretely, in the fourth mesh level, we obtain $25\%$ of the elements with an energy error associated with order $p-1$, or order $p-2$. The error of these elements is much higher than the elements with an energy error related to $p$ order approximation, making the error of these $p-1$ or $p-2$ elements dominant compared to the rest of the elements.

Given the results of these experiments, we conclude that the non-uniform $p$ refinement implementation is correct. Then, in the following experiments, we focus on completing the analysis of our implementation to solve more challenging problems.

\subsection{$h-$adapted meshes for a synthetic marine layered model}
\label{test_layered_meshing_rules}
As a second example, we investigate the performance of the $h-$adapted meshes. More concretely, we perform a convergence test to determine the number of points per skin depth $\lambda_{\delta}$ and the resolution number $r_{s}$ which are required to design $h-$adapted meshes that satisfy a given error threshold. For this experiment, we consider an VTI synthetic layered model that consists of $3.5\,$km of air layer ($\rho_{\text{air}}=\num{e8}\,\Omega \cdot\,$m), 1$\,$km of seawater ($\rho_{\text{sea}}=0.3\, \Omega \cdot$m), followed by 2$\,$km of sediments with VTI ($\rho_{\text{h}}=4\,\Omega \cdot$m and $\rho_{\text{v}}=6\,\Omega \cdot$m), and finally 1$\,$km of isotropic basement ($\rho_{\text{basement}}=1000\,\Omega \cdot$m). We execute our simulations in parallel with 480 processors. The system of equations is solved using the \texttt{PETSc}\xspace implementation of the \textit{GMRES} solver. The reference solution for a 3$\,$Hz $x$-directed dipole transmitter is computed semi-analytically using the \empymod code~\citep{GEO.17.Werthmuller}.

\begin{figure*}
  \centering
  \begin{subfigure}[b]{0.45\textwidth}
    \includegraphics[trim= 0 0 0 0, clip,width=\textwidth]{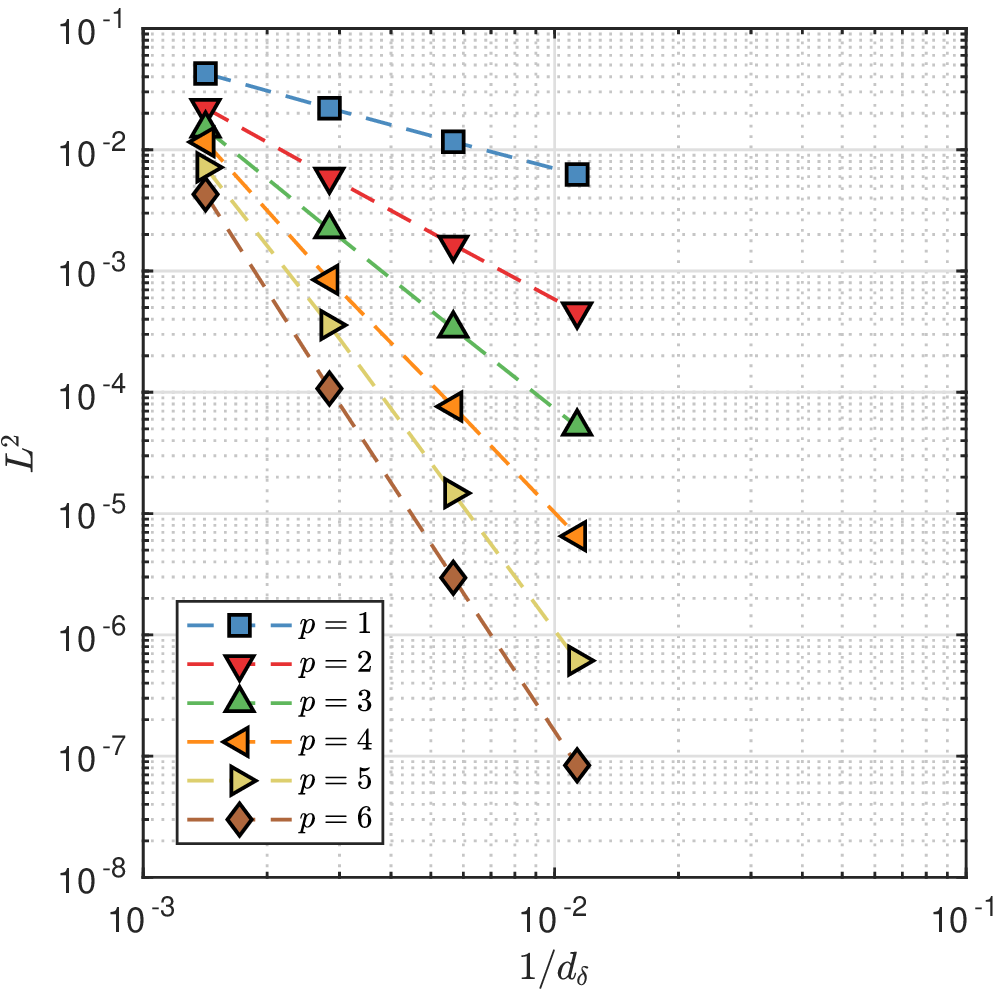}
    \caption{}
    \label{fig:convergence_a}
  \end{subfigure}
  \begin{subfigure}[b]{0.45\textwidth}
    \includegraphics[trim= 0 0 0 0, clip, width=\textwidth]{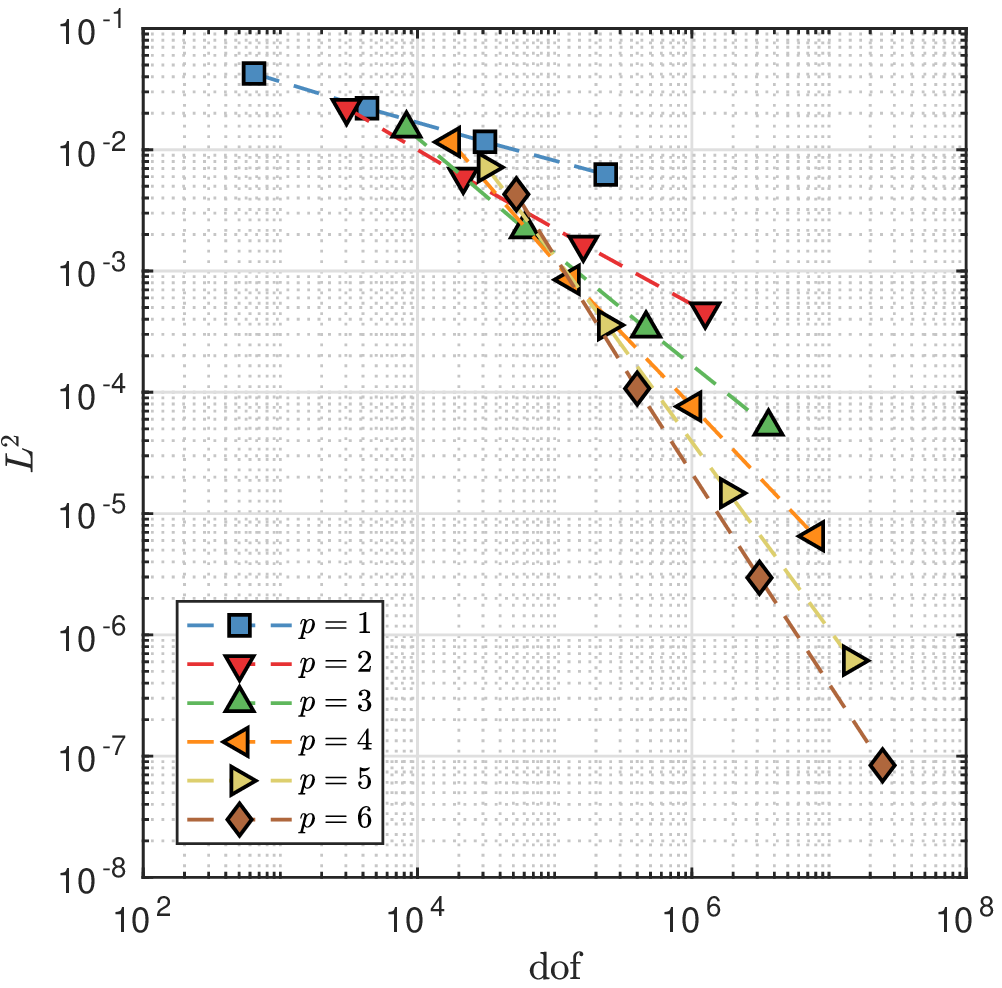}
    \caption{}
    \label{fig:convergence_b}
  \end{subfigure}
  \begin{subfigure}[b]{0.45\textwidth}
    \includegraphics[trim= 0 0 0 0, clip, width=\textwidth]{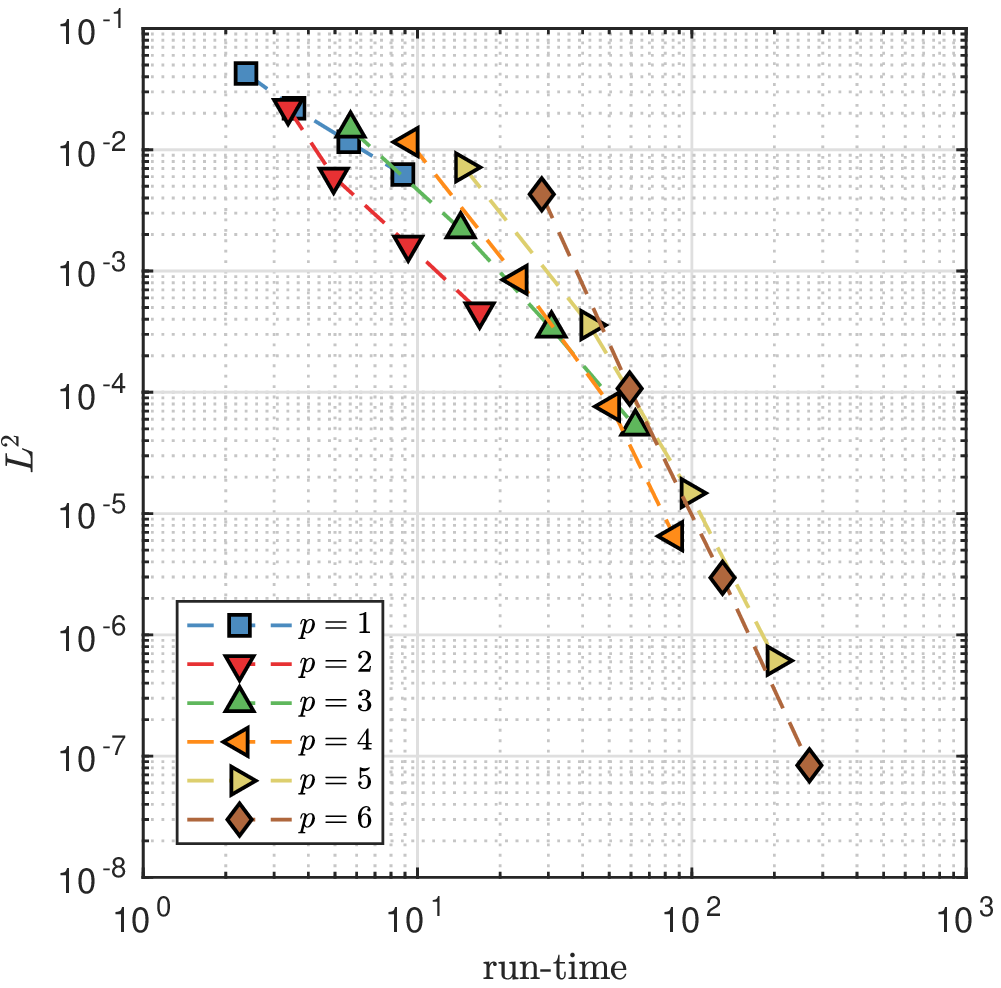}
    \caption{}
    \label{fig:convergence_c}
  \end{subfigure}
  \begin{subfigure}[b]{0.45\textwidth}
    \includegraphics[trim= 0 0 0 0, clip, width=\textwidth]{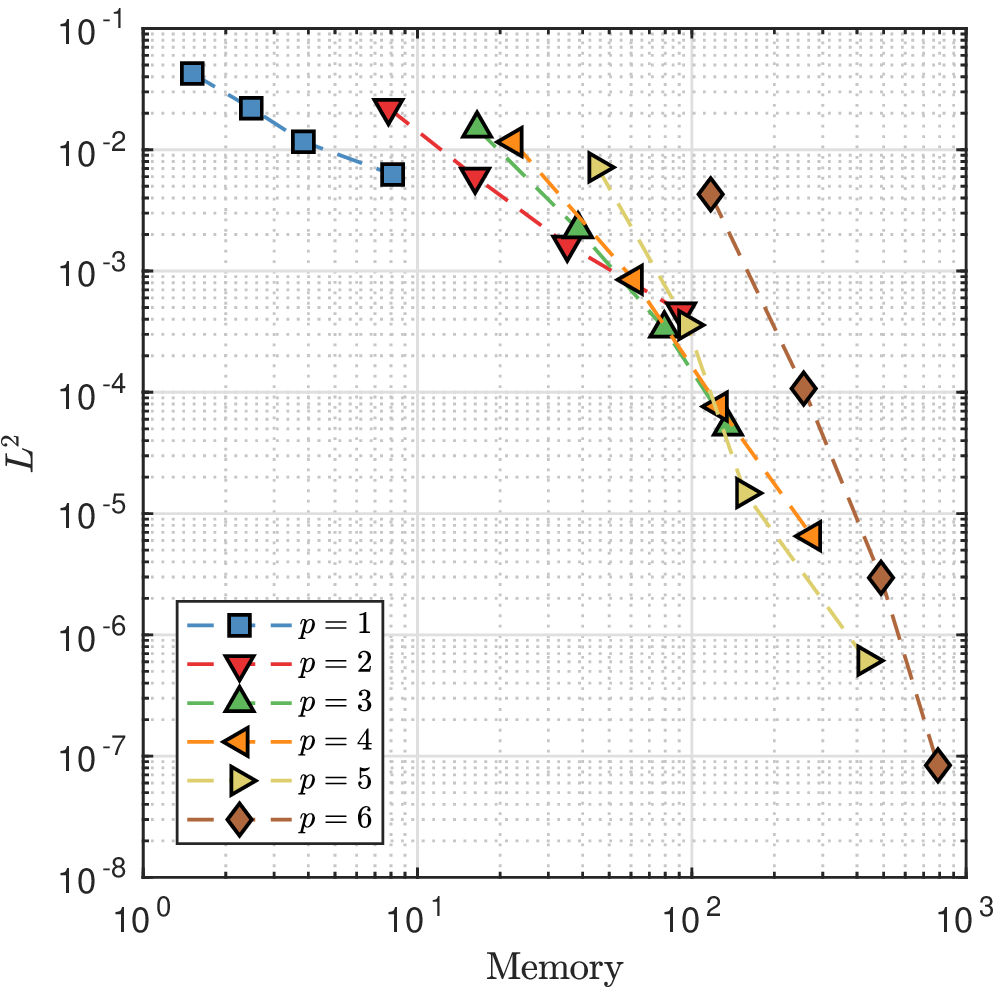}
    \caption{}
    \label{fig:convergence_d}
  \end{subfigure}
  \caption{Convergence behavior in terms of  (a) characteristic mesh spacing $d_\delta$, (b) number of dof, (c) run-time, and (d) memory consumption.}
\end{figure*}

We design a set of meshes with global hierarchical $h-$refinement, whose resulting statistics are given in Table~\ref{table:mesh_hierarchies_layered}. Then, on these meshes, we compute the 3D CSEM solutions for each polynomial basis degree (e.g., $p=1,2,3,4,5,6$). The obtained $L^{2}$ errors, convergence rates, run-time, and memory consumption are summarized in Table~\ref{table:convergence_layered}. For all nested refined meshes, the basis order $p=6$ produces the most accurate solution, the approximations $p=5, 4,3,2$ follows, and the piecewise $p=1$ is the least accurate solution, as expected. Regarding mesh convergence, Fig.~\ref{fig:convergence_a} shows the rates depicted in Table~\ref{table:convergence_layered}, where it can be seen that convergence slopes are consistent with each polynomial basis order. Also, Fig.~\ref{fig:convergence_b} depicts the trade-off between numerical error and dof. Despite number of dof grows faster for high-order polynomial basis, this increment is compensated by the decrease in the number of tetrahedral elements required to reach a given numerical error level (e.g., for $p=2$ the error obtained with the second mesh level is almost similar than $p=3$ with same mesh level, and much better than $p=1$ with third mesh). However, a better pattern is observed for low-order transitions (from $p=1$ to $p=2,3$) than for high-order transitions (from $p=3$ to $p=4,5,6$). Still, in function of the desirable accuracy, the basis orders $p=3,4$ can be competitive (e.g., for $p=3$ the obtained error with the fourth mesh level is similar than those achieved with $p=4$ and the third mesh level). Furthermore, a close inspection of Fig.~\ref{fig:convergence_c} shows that high-order piecewise demands more run-time to reach the final solution for the same mesh level. Finally, Fig.~\ref{fig:convergence_d} shows the memory consumption for each polynomial approximation, where it can be seen that high-order basis demands more memory to solve the problem under consideration. 

Given each approximation scheme's reasonable performance ranges (e.g., Table~\ref{table:convergence_layered}), we conclude that, when not using non-uniform $p$ refinement, the $p=2$ piecewise provides the best general trade-off between mesh size, accuracy, and computational effort. However, the optimum choice of the polynomial scheme depends on the input model. This conclusion is consistent with those reported in previous studies~\citep{Grayver2015,Schwarzbach2011,GJI.19.CastilloReyes}.

Based on our convergence results, we obtain the number of points per skin depth $\lambda_{\delta}$ and resolution number $r_{s}$ required to design $h-$adapted meshes that satisfy 5, 3, and 1 percent error thresholds. These mesh resolution parameters are depicted in Table~\ref{table:points_skindepth}. Using these meshing rules, the resulting global characteristic mesh size and local $h-$refinement ensure the proper resolution of the EM wave in the computational domain. Last but not least, we acknowledge that although $\lambda_{\delta}$ and $r_{s}$ values are suitable for designing unstructured adapted meshes of sufficient discretization quality, the best performance compromise depends on the input model (similar conclusion to that presented by~\cite{Plessix2007}).
\begin{table}
\caption{Mesh statistics for a convergence test on a marine CSEM model with VTI. Hierarchical mesh level, number of elements, characteristic mesh spacing $d_{\delta}$ (expressed in meters), and number of dof for each HEFEM polynomial degree $p=1,2,3,4,5,6$.}
\begin{center}
\footnotesize
\begin{tabular}{c c c c c c c c c}
\toprule
\multirow{2}{*}{Mesh level} & \multirow{2}{*}{Elements} & \multirow{2}{*}{$d_{\delta}$} & \multicolumn{6}{c}{dof} \\ \cmidrule{4-9}
        & &  & $p=1$ &  $p=2$ & $p=3$ & $p=4$ & $p=5$ & $p=6$\\
\midrule
1 &  $373$      & $705.578$  & $643$ & $3\,038$ & $8\,304$ & $17\,560$ & $31\,925$ & $52\,518$ \\
2 &  $2\,984$   & $352.789$  & $4\,287$ & $21\,550$ & $60\,741$ & $130\,812$ & $240\,715$ & $399\,402$\\
3 &  $23\,872$  & $176.394$  & $31\,022$ & $161\,692$ & $463\,626$ & $1\,008\,440$ & $1\,867\,750$ & $3\,113\,172$ \\
4 &  $190\,976$ & $88.197$   & $235\,388$ & $1\,251\,320$ & $3\,620\,724$ & $7\,916\,528$ & $14\,711\,660$ & $24\,579\,048$\\
\bottomrule
\end{tabular}
\end{center}
\label{table:mesh_hierarchies_layered}
\end{table}
\begin{table}
\caption{Convergence results for each polynomial degree $p=1,2,3,4,5,6$ over each mesh level depicted in Table~\ref{table:mesh_hierarchies_layered}. Numerical solution error in $L^{2}$-norm, convergence orders $\mathcal{O}_{L^{2}}$, run-time (minutes), and memory consumption (Gb) are included.}
\begin{center}
\begin{tabular}{c c c c }
\toprule
$L^{2}$  &  $\mathcal{O}_{L^{2}}$ & 
Run-time & Memory\\
\midrule 
\multicolumn{4}{c}{$p=1$} \\
\num{4.252e-2} &   $-$  & 2.369 & 1.510 \\
\num{2.201e-2} &  0.949 & 3.548 & 2.478 \\ 
\num{1.161e-2} &  0.918 & 5.615 & 3.835 \\
\num{6.301e-3} &  0.897 & 8.841 & 8.115 \\

\multicolumn{4}{c}{$p=2$} \\
\num{2.210e-2}  &   $-$   &  3.374 &  7.832 \\
\num{6.010e-3}  &  1.882  &  4.944 & 16.216 \\
\num{1.601e-3}  &  1.859  &  9.247 & 35.147 \\
\num{4.630e-4}  &  1.833  & 16.841 & 91.310 \\

\multicolumn{4}{c}{$p=3$} \\
\num{1.495e-2} &    $-$ &   5.698 &  16.487 \\
\num{2.231e-3} &  2.756 &  14.364 &  38.478 \\
\num{3.359e-4} &  2.718 &  30.781 &  79.364 \\ 
\num{5.216e-5} &  2.687 &  62.247 & 135.369 \\

\multicolumn{4}{c}{$p=4$} \\
\num{1.161e-2}  &    $-$  &  9.354  &  22.364 \\
\num{8.470e-4}  &   3.772 &  23.365 &  61.391 \\ 
\num{7.622e-5} &   3.473 &  50.784 & 125.369 \\ 
\num{6.501e-6} &   3.551 &  86.147 & 274.144 \\

\multicolumn{4}{c}{$p=5$} \\
\num{7.201e-3}  &     $-$  &  14.784  &  45.321 \\
\num{3.582e-04} &   4.312 &  42.369 &  96.259 \\
\num{1.474e-05} &   4.603 &  98.178 &  156.367 \\
\num{6.125e-07} &   4.588 &  202.378 &  433.412 \\

\multicolumn{4}{c}{$p=6$} \\
\num{4.312e-3}  &     $-$  &  28.368  &  117.140 \\
\num{1.073e-04} &   5.320 &  59.327 &  255.690 \\
\num{2.956e-06} &   5.182 &  129.382 & 489.412 \\
\num{8.396e-08} &   5.137 &  268.328 & 789.291 \\

\bottomrule
\end{tabular}
\end{center}
\label{table:convergence_layered}
\end{table}
\begin{table}
\caption{Number of points per skin depth $\lambda_{\delta}$ and resolution number $r_{s}$ required to design $h-$adapted meshes that satisfy 5, 3, and 1 percent error thresholds. These parameters are given for each polynomial order $p$.}
\begin{center}
\begin{tabular}{c c c c c c c c c}
\toprule
\multirow{2}{*}{$p$} & \multicolumn{2}{c}{$5\%$} & \phantom{abc} & \multicolumn{2}{c}{$3\%$} & \phantom{abc} & \multicolumn{2}{c}{$1\%$}\\ \cmidrule{2-3} \cmidrule{5-6} \cmidrule{8-9}
 &  $\lambda_{\delta}$ & $r_{s}$ & \phantom{abc} & $\lambda_{\delta}$ & $r_{s}$ & \phantom{abc} & $\lambda_{\delta}$ & $r_{s}$\\
\midrule 
$1$  & 9.22 & 14 & & 13.93 & 13 & & 16.81 & 15\\
$2$  & 7.12 & 10 & & 9.39 & 11 & & 10.55 & 12 \\
$3$  & 6.33 &  9 & & 8.14 & 10 & & 9.35 & 11 \\
$4$  & 5.02 &  8 & & 5.95 &  9 & & 8.12 & 10\\
$5$  & 4.98 &  7 & & 5.82 &  8 & & 6.95 & 9 \\
$6$  & 3.86 &  6 & & 4.69 &  7 & & 5.78 & 8 \\
\bottomrule
\end{tabular}
\end{center}
\label{table:points_skindepth}
\end{table}

\subsection{$p-$refinement experiment for the marine MR3D model}
\label{marine_mr3d_model}
\begin{figure}
	\centering
	\includegraphics[trim= 0 0 100 80,clip,width=0.55\textwidth]{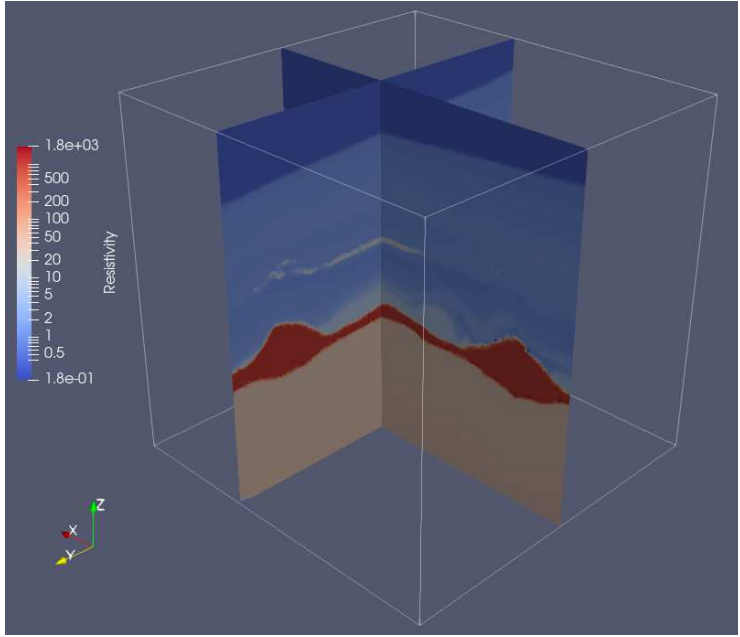}
	\caption{Marine MR3D resistivity model proposed and released under the open creative common license (CC by 4.0) by~\cite{GEO.19.Correa}.}
	\label{fig:fig3}
\end{figure}

The third example focuses on studying the impact of uniform $p-$refinement on a realistic 3D resistivity model with VTI. The model under consideration was proposed and released under the open creative common license (CC by 4.0) by~\cite{GEO.19.Correa}. The available model data, also referred to as Marlim R3D (MR3D), consists of EM responses for six frequencies from $0.125\,$Hz to $1.25\,$Hz. The up-scaled model compromises $515 \times 563 \times 310$ cells where each cell has dimensions of $100 \times 100 \times 20\,$m, resulting in almost $90$ million cells. The data-acquisition area is a regular grid composed of 500 stations placed $1\,$km above the irregular seafloor. Although the original published data includes EM responses for 45 transmitter-towlines, we only compute the electric field for the in-line transmitter \texttt{04Tx013a} and the cross-line transmitter \texttt{04Tx014a} on the receivers profile \texttt{04Rx251a}. For both transmitters, the frequency of the excitation current is $1.25\,$Hz. Fig.~\ref{fig:fig3} shows a 3D view of the conductivity model under consideration. We used the FE EM responses provided by~\cite{GJI.21.Werthmuller} as reference solution. For more details about the oil field MR3D model, we refer to~\cite{GEO.19.Correa}.

For this marine 3D CSEM modeling, we investigate the impact of global and uniform $p-$refinement on EM responses. The experiments are performed in two phases. First, we design $h-$adapted meshes for $p=1,2,3,4,5,6$. Second, for each polynomial scheme, we compute EM responses applying global and uniform $p-$refinement at different mesh entity levels:
\begin{enumerate}
    \item Edge level: only on edges.
    \item Face level: only on faces.
    \item Interior, or volume level: inside each element.
    \item Element level: in all mesh element entities (e.g., edges, faces, and volume per element). 
\end{enumerate}
The resulting mesh statistics for each piecewise and each refinement entity are given in Table~\ref{table:marlim_mesh_statistics}. Each numerical approximation has been computed in parallel with 480 cores and using the direct \texttt{MUMPS}\xspace solver.
\begin{table}
\scriptsize
\caption{$h-$adapted mesh statistics for MR3D model. Piecewise $p$, number of elements, and number of dof for $p-$refinement at different mesh entity level (edge, face, volume, and element). Run-time (minutes) and memory consumption (Gb) are also given.}
\begin{center}
\begin{tabular}{c c c c c c c c c c}
\toprule
\multirow{2}{*}{$p$} & \multirow{2}{*}{Elements} & \multicolumn{2}{c}{Edge} & \multicolumn{2}{c}{Face}  & \multicolumn{2}{c}{Volume} & \multicolumn{2}{c}{Element}   
\\ 
\cmidrule(l){3-4} \cmidrule(l){5-6} \cmidrule(l){7-8} \cmidrule(l){9-10}\\
& & dof & Run-time/Memory & dof & Run-time/Memory & dof & Run-time/Memory & dof & Run-time/Memory \\
\midrule 
1 & $14\,430\,811$ & $15\,241\,416$ & 525.44/572.45 & $-$ & $-$ & $-$ & $-$ &  $15\,241\,416$ & 572.5/525.1\\
2 & $916\,992$     & $2\,199\,118$  & 75.81/82.59 & $4\,815\,901$ & 180.88/166.02   & $-$ & $-$ &  $5\,915\,460$  & 222.2/203.9\\
3 & $666\,726$ & $2\,408\,487$ & 83.03/90.46 & $8\,922\,323$  & 335.11/307.59 &  $2\,803\,007$ & 105.27/96.63 & $12\,528\,159$ & 407.5/431.9\\
4 & $464\,483$ & $2\,246\,756$ & 77.45/84.38 & $11\,893\,097$ & 446.69/410.01 & $6\,135\,485$  & 230.44/211.52 & $19\,151\,960$ & 719.3/660.3\\
5 & $311\,030$ & $1\,893\,170$ & 65.26/71.10 & $13\,057\,734$ & 490.43/450.16 & $9\,709\,534$  & 364.68/334.73 & $23\,903\,170$ & 897.8/824.1\\
6 & $248\,243$ & $1\,818\,486$ & 62.69/68.30 & $15\,501\,231$ & 582.21/534.40 & $15\,197\,661$ & 570.81/523.93 & $31\,911\,216$ & $1\,198.6$/$1\,100.1$ \\
\bottomrule
\end{tabular}
\end{center}
\label{table:marlim_mesh_statistics}
\end{table}

\begin{figure}
	\centering
	\begin{subfigure}[b]{0.45\textwidth}
        \includegraphics[trim= 0 0 0 0,clip,width=\textwidth]{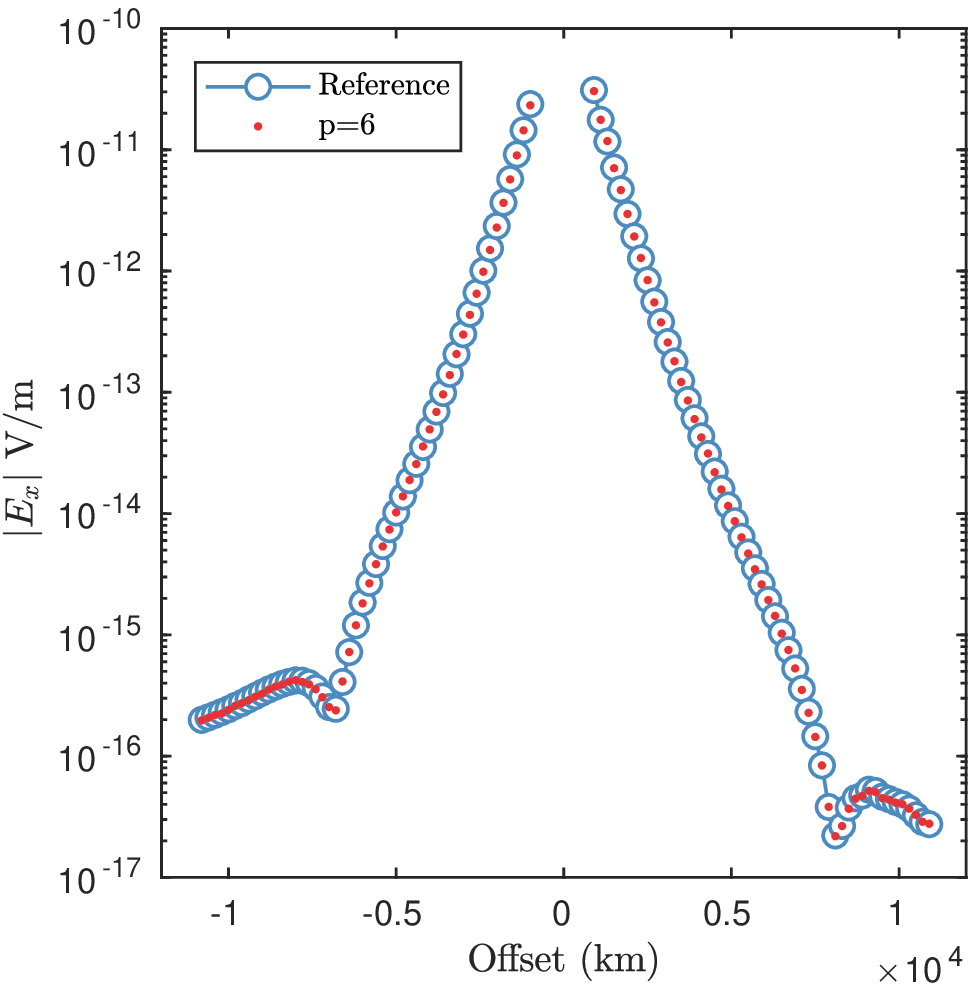}
        \caption{}
        \label{fig:fig4a}
    \end{subfigure}
    \hspace{.2cm}
    \begin{subfigure}[b]{0.45\textwidth}
        \includegraphics[trim= 0 0 0 0,clip,width=\textwidth]{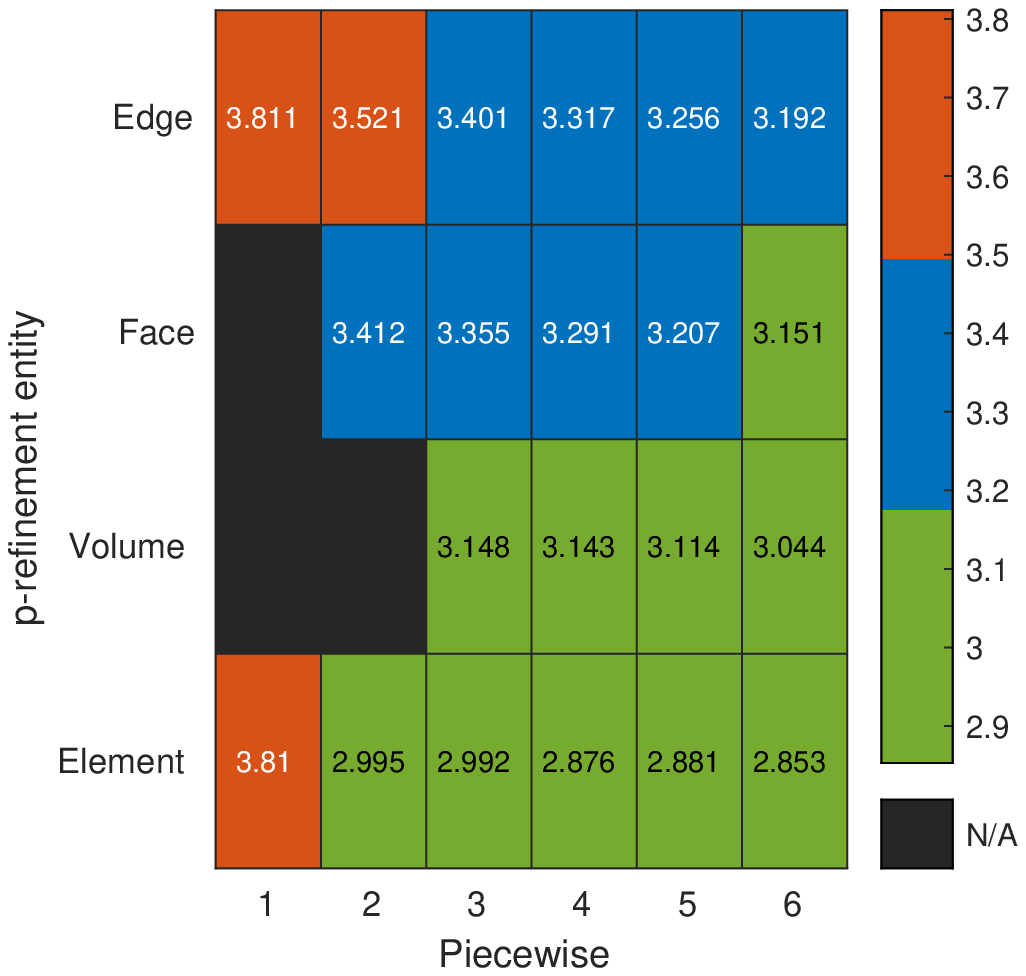}
        \caption{}
        \label{fig:fig4b}
    \end{subfigure}
    \caption{Comparison of $x-$component amplitude between \petgem and reference data. (a) EM responses for in-line transmitter \texttt{04Tx013a} and basis order $p=6$ (similar results have been obtained for $p=1,2,3,4,5$). (b) Maximum misfits percentile for each polynomial order and each uniform $p-$refinement entity (edge, face, volume, element). The patch is organized piecewise in x-axis and refinement level in y-axis. Black colors indicated non applicable (N/A) refinement in that entity.}
	\label{fig:fig4}
\end{figure}

As aforementioned, we computed synthetic EM fields for each piecewise. Still, to preserve brevity, Fig.~\ref{fig:fig4a} only shows $p=6$ EM responses for the in-line transmitter \texttt{04Tx013a} and using uniform refinement at element level (e.g., dof in edges, faces, and volume of each tetrahedral element). Here, the EM amplitudes are almost identical with respect to the reference. Similar results have been obtained for $p=1,2,3,4,5$ at element level refinement. Fig.~\ref{fig:fig4b} depicts the maximum misfits percentile for each polynomial order and each uniform $p-$refinement entity, which are close to the error threshold used to design each mesh ($3\%$ percent error). Regardless of the refined entity, the basis order $p=6$ produces the most accurate EM responses, followed by the $p=5,4,3,2,1$ approximations. More concretely, if uniform refinement is only applied on edges, the maximum error is reduced from $3.811\%$ for $p=1$ to $3.192\%$ for $p=6$. Similar behaviour is observed for uniform refinement on faces where the error is decreased from $3.412\%$ for $p=2$ to $3.151\%$ for $p=6$. For uniform refinement on volume (e.g., interior dof), the error is reduced $3.148\%$ for $p=3$ to $3.044\%$ for $p=6$. However, a close inspection of misfits depicted in Fig.~\ref{fig:fig4b} shows that uniform refinement on volume produces more accurate solutions than on edges or faces. Furthermore, the obtained misfits for uniform $p-$refinement on volume is close to those obtained at for uniform $p-$refinement at element level (e.g., for $p=3$, the obtained misfit at volume level is close to that obtained at element level ($3.148\% \approx 2.992\%$)). Last but not least, we acknowledge that uniform $p-$refinement at element level offers the best misfit ratios. 

\begin{figure*}
	\centering
	\begin{subfigure}[b]{0.45\textwidth}
        \includegraphics[trim= 0 0 0 0,clip,width=\textwidth]{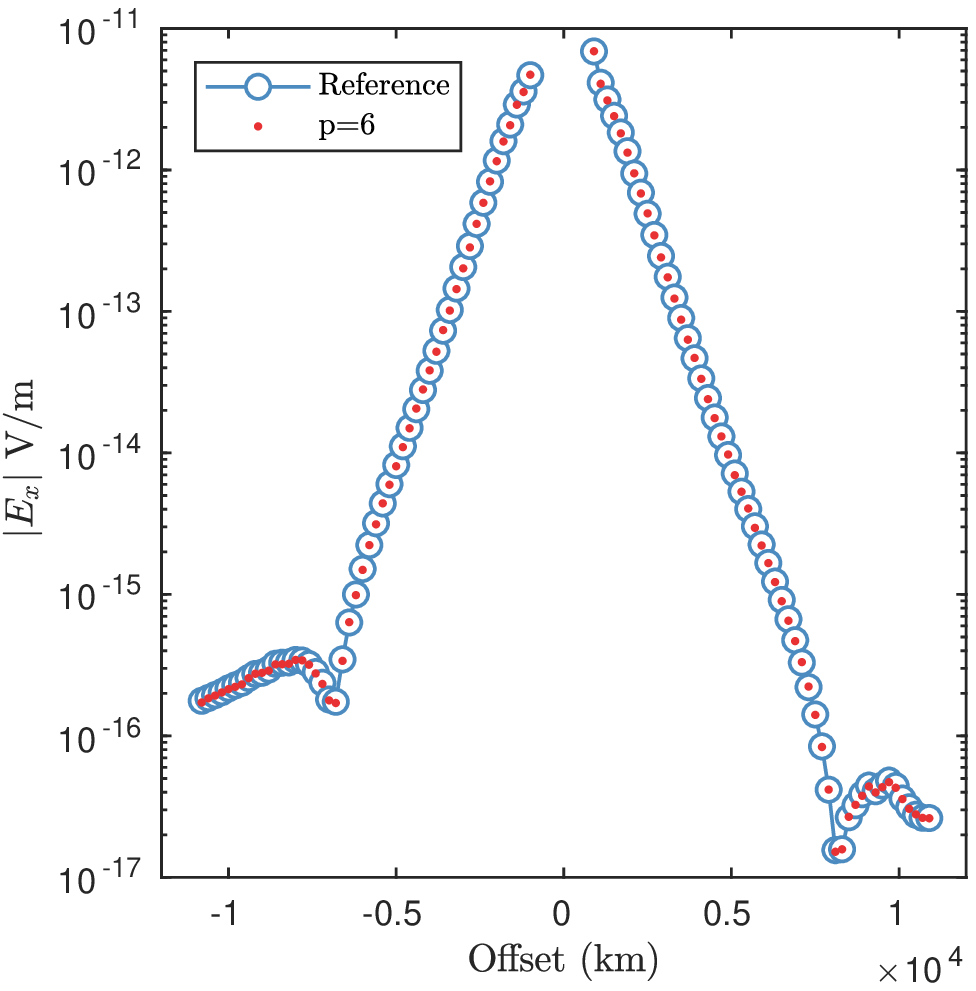}
        \caption{}
        \label{fig:fig5a}
    \end{subfigure}
    \hspace{.2cm}
    \begin{subfigure}[b]{0.45\textwidth}
        \includegraphics[trim= 0 0 0 0,clip,width=\textwidth]{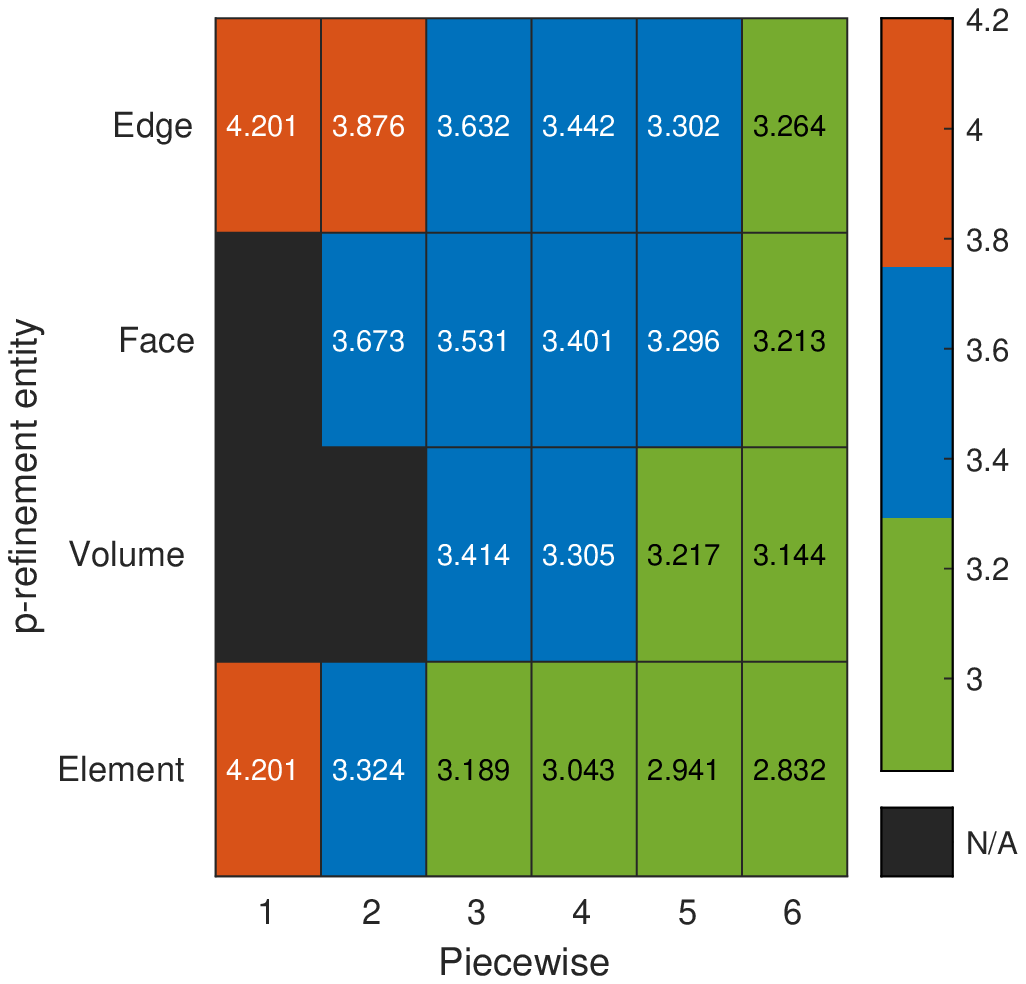}
        \caption{}
        \label{fig:fig5b}
    \end{subfigure}
    \caption{Comparison of $x-$component amplitude between \petgem and reference data. (a) EM responses for in-line transmitter \texttt{04Tx014a} and basis order $p=6$ (similar results have been obtained for $p=1,2,3,4,5$). (b) Maximum misfits percentile for each polynomial order and each uniform $p-$refinement entity (edge, face, volume, element). The patch is organized piecewise in x-axis and refinement level in y-axis. Black colors indicated non applicable (N/A) refinement in that entity.}
	\label{fig:fig5}
\end{figure*}

Fig.~\ref{fig:fig5a} shows the EM responses for the cross-line transmitter \texttt{04Tx014a} computed with $p=6$. Fig.~\ref{fig:fig5b} depicts the maximum misfit percentile for each piecewise and each uniform $p-$refinement entity. The conclusions for this cross-line profile are similar to those obtained for previous in-line profile analysis.

Given our numerical results, we conclude that dof on volume level has the most significant impact on reducing the error. The dof on volume are followed by dof on faces, and finally by dof on edges. Furthermore, although a slight improvement in the error is observed during the transition from $p=2$ to $p=3$, this represents a considerable increment in computational cost to obtain the solution (e.g., the number of dof increases from $5\,915\,460$ to $12\,528\,159$, but the error is practically the same $2.995 \approx 2.992$). This pattern is also preserved for higher polynomial basis. Then, this analysis confirms that second-order uniform $p-$refinement in conjunction with $h-$adapted meshes offer the best compromise between accuracy and run-time. We point out that there might be other $p-$refinement strategies combination (e.g., based on guided error estimators) that can solve the 3D CSEM modeling more efficiently. However, these schemes are beyond our paper's scope as we focus on analyzing whether the use of HEFEM with $h-$adapted meshes and supervised $p-$refinement strategies is justified.
\subsection{$p-$refinement experiment for the land Vall\'es model}
\label{land_valles_model}
\begin{figure}
	\centering
	\includegraphics[trim= 0 0 0 0,clip,width=0.80\textwidth]{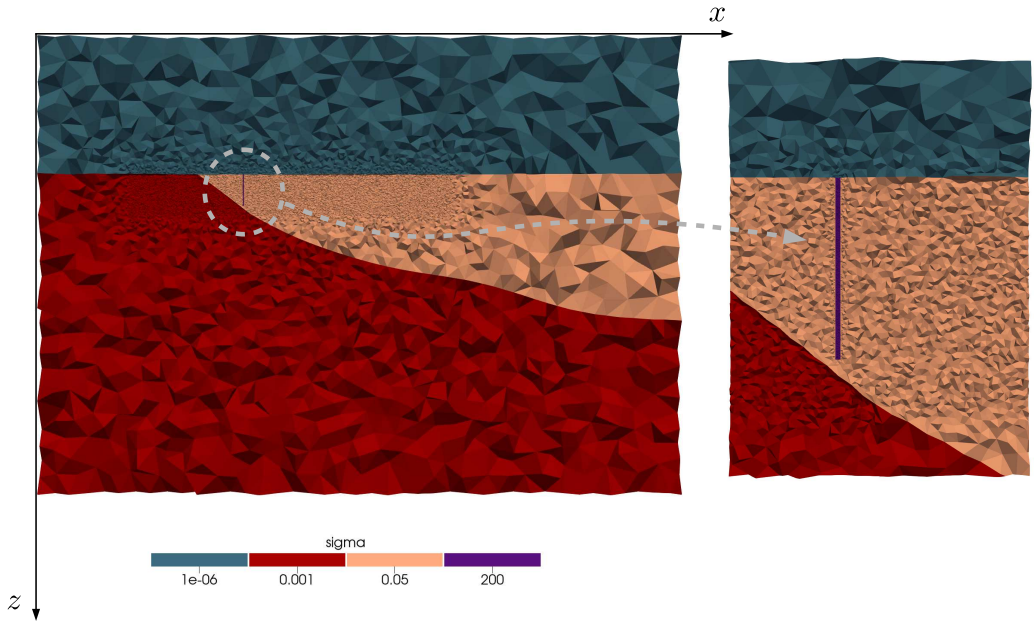}
	
	\caption{3D CSEM Vall\'es model with its resulting $h-$adapted unstructured tetrahedral mesh. The 2D view corresponds to $y=2$km with zoom-in closer to metallic casing vicinity.}
	\label{fig:fig6}
\end{figure}

\begin{figure}
	\centering
	\includegraphics[trim= 0 0 0 0,clip,width=0.80\textwidth]{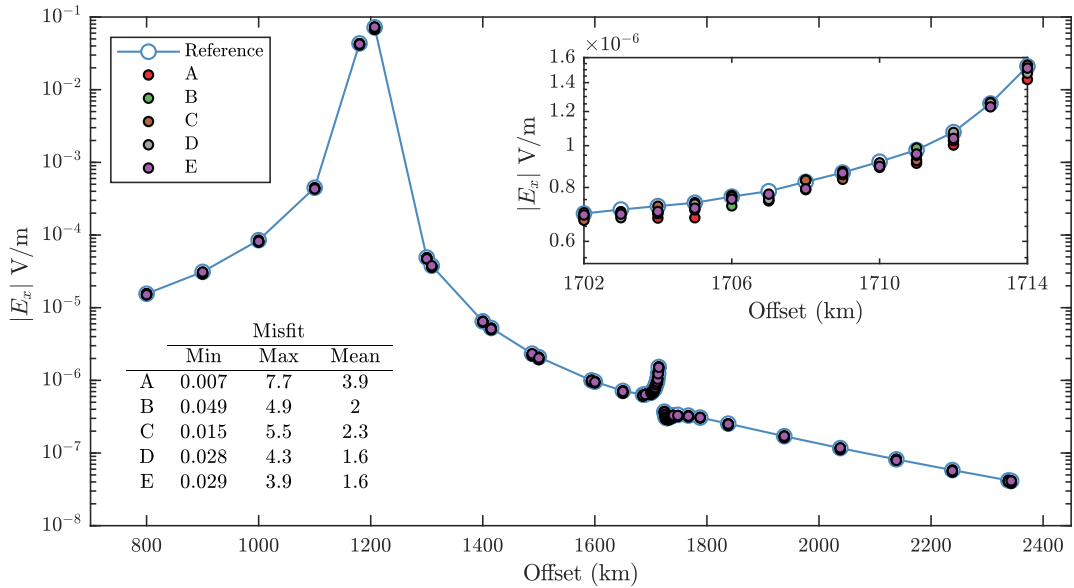}
	
	\caption{EM responses for the land Vall\'es model using non-homogeneous $p-$refinement. The misfits for each refinement scheme are also included.}
	\label{fig:fig7}
\end{figure}

Finally,  the  fourth  test  focuses  on  studying  the  impact  of non-uniform $p-$refinement on  the  EM  responses.  To perform these experiments we use a 3D land CSEM model in the geothermal exploration context. The model under consideration, referred to as Vall\'es model and introduced by~\cite{TGRS.21.CastilloReyes}, is composed by 4$\,$km thick air layer ($\rho_{\text{air}}=\num{e8}\,\Omega \cdot\,$m), resistive basement ($\rho_{\text{basement}} = 1000\,\Omega \cdot\,$m) including topography, and the conductive sediments ($\rho_{\text{sediments}} = 20\,\Omega \cdot\,$m). To simulate the metallic old boreholes present in the area, the model includes a vertical cylinder ($\rho_{\text{casing}} = \num{e-4}\,\Omega \cdot\,$m) embedded in the sediments. The cylinder length is $200\,$m, and it is centered at $x=1.719\,$km, $y=2\,$km, and $z=-100\,$m. The Vall\'es model is depicted in Fig.~\ref{fig:fig6}. We use an $x$-directed transmitter
with a moment of $1\,$Am and frequency of $2\,$Hz. This source is located at $x=1.209\,$km, $y=2\,$km, and $z=0\,$km. The data-acquisition consists of eighty-one stations placed in-line to the transmitter position and along its orientation. The computed responses for a metallic casing distant to the $x-$directed transmitter are shown in the original publication and compared against experimental data in-situ. We reproduce the EM responses for this set-up in our experiments. The resulting systems of equations are solved in parallel with 480 cores and using the \texttt{MUMPS}\xspace solver. 
For this model, we simulate different non-uniform $p-$refinement configurations that allow us a detailed analysis of the most relevant cases. Table~\ref{table:valles_mesh_statistics} shows the experiment details. Fig.~\ref{fig:fig7} compare the electric field amplitude of the $x-$component for each $p-$refinement scheme. A close inspection of these EM field responses confirms the effect of the metallic casing presence, which is most noticeable close to the borehole vicinity (from $x=1\,700\,$m to $x=1\,740\,$m). This behavior is observed regardless of the $p-$refining method used. Fig.~\ref{fig:fig7} also depicted the obtained misfits for each test. Here, the cases D and E produce the most accurate solutions ($1.6\%$ of average misfit). However, the solution E is not competitive if its number of dof is compared to the strategy D ($\approx 7$ millions of dofs for case E against $\approx 3$ millions of dofs for case D). The solution B is the third most accurate. The least accurate solutions are the C and A options. Given these results, we conclude that when it is desired to use non-uniform $p-$refinement, the dof associated with the element volume have a greater impact on controlling the solution's error (e.g., case D implemented $p=1$ in edges and faces, and $p=3$ in volume). However, for uniform $p-$refinement, the $p=2$ offers the best compromise between accuracy and computational effort. This conclusion is consistent with previously published results~\citep{Schwarzbach2011,Grayver2015}. Furthermore, our numerical results suggest that a good compromise between $h-$adapted meshes and $p-$refinement is needed to provide accurate solutions in a feasible run-time.
\begin{table}
\caption{Summary of scenarios studied for Vall\'es model. The dash means that there are no dof associated with that entity (equivalent to $p=1$). Run-time (minutes) and memory consumption (Gb) are also given.}
\begin{center}
\begin{tabular}{c c c c c c}
\toprule
\multirow{2}{*}{Case label} & \multicolumn{3}{c}{$p-$refinement} & \multirow{2}{*}{dof} & \multirow{2}{*}{Run-time/Memory}\\
\cmidrule{2-4}
 &  Edge & Face & Volume & \\
\midrule 
A & $p=1$ & $-$ & $-$ & $939\,258$ & 32.27/32.38\\
B & $p=2$ & $p=2$ & $-$ & $5\,094\,536$ & 191.34/175.63\\
C & $p=1$ & $p=2$ & $-$ & $4\,155\,278$ & 156.06/143.25\\
D & $p=1$ & $p=1$ & $p=3$ & $3\,341\,784$ & 125.51/115.20\\
E & $p=2$ & $p=2$ & $p=3$ & $7\,497\,062$ & 281.58/258.46\\
\bottomrule
\end{tabular}
\end{center}
\label{table:valles_mesh_statistics}
\end{table}

\section{Conclusions}
\label{conclusions}
The development of geo-electromagnetic modeling routines has increased in the last years, and today there are several options to obtain reasonable-looking results. However, with the sole exception of the algorithms developed by~\cite{Plessix2007,Schwarzbach2011,Grayver2015,GJI.19.CastilloReyes,Rochlitz2019} and regardless of the type of meshes used, most of these 3D modeling routines do not employ tailored meshes and do not support high-order discretizations. These are the core motivations for this study: to evaluate the benefits and limitations of HEFEM and supervised $h+p$ tailored meshing to solve challenging 3D CSEM setups.

We perform numerical experiments as a proof-of-principle for the computation of EM responses based upon HEFEM, supervised $h+p$ refinement strategies, and HPC. These  experiments provide relevant information about the performance of high-order polynomials ($p=1,2,3,4,5,6$) in conjunction with unstructured and tailored meshes. As far as the authors know, only $p=1,2,3$ performance metrics are reported in previous publications. Furthermore, we consider 3D CSEM models with synthetic and experimental data, making them valid for academic and industrial applications. The points mentioned above are clear contributions with respect to the rest of the state-of-the-art works. 

To demonstrate the robustness of our numerical schemes, we compute the solutions for marine and land test cases with real-life applications (e.g., mineral mining, oil $\&$ gas, and geothermal reservoir characterization and interpretation). We focus on studying its compromise in terms of convergence rate, run-time, and memory needs. In our numerical experiments, the accuracy obtained with each polynomial function is consistent with the theoretical definition. Also, the high-order polynomial degrees require a fewer number of dof to attain a given error level in comparison with the low-order case. However, this accuracy improvement has a cost. The computational implementation of HEFEM and its parallelization is technically complex. Also, the use of high-order elements decreases the sparsity pattern of the resulting linear system due to the number of dof per element is larger~\cite{Schwarzbach2011,Grayver2015}. Then, high-order schemes demand more memory. Furthermore, the computation of the element integrals is more expensive for high-order piecewise (e.g., the quadrature order for numerical integration increases in a proportional ratio to the polynomial basis order). Consequently, the run-time for linear system assembly is also increased. However, the linear solver is still the most demanding phase in run-time and memory consumption. Therefore, it is necessary to have efficient assembly routines and robust and scalable linear equation solvers to face these challenges. We mitigated the challenges mentioned above by using an HPC code with general-purpose solvers. 

Regarding $h-$refinement, we have implemented a supervised strategy to design $h-$adapted meshes for each basis degree. These meshing rules control the global characteristic mesh size and local $h-$refinement to ensure a proper discretization of the EM problem under consideration. Accurate synthetic EM responses were obtained by using this meshing strategy in conjunction with high-order basis. Regarding $p-$refinement, our HEFEM modeling routine supports arbitrary polynomial degrees at each element (e.g., edge, face, volume). From an academic and theoretical perspective, these high-order schemes are usually favorable. However, they are more expensive to deal with from a practical perspective. Based on the results presented here, uniform refinement for high-order polynomial degrees ($p>3$) are more resource-demanding (e.g., run-time and memory), making them less attractive or limiting for modest computational architectures. A $p=2$ uniform refinement exhibits the best accuracy/performance trade-off for the modeling test cases studied here. This conclusion is consistent with those reported by~\cite{Farquharson2011} and~\cite{Grayver2015}. However, its performance compromise depends on the input model (e.g., frequency, conductivity/resistivity contrasts, solver type, mesh quality, domain decomposition strategy, among others). For non-uniform refinement, the dof associated to the element volume has the most significant impact on controlling the misfits on EM responses. The dof associated to element faces and element edges follows, respectively. In our numerical results, the refinement scheme with $p=3$ on volume and $p=1$ on faces and edges has exhibited the best compromise between accuracy and computing effort. However, its misfits are comparable to those obtained with a uniform refinement using a second-order approximation (e.g., see Fig.~\ref{fig:fig6}). Then, we conclude that a non-uniform $p-$refinement scheme can be competitive if applied hierarchically. First, increase the $p$ order in volume, then in faces, and finally in edges. Then, we can select preferable modelling options with our analyses:
\begin{enumerate}
    \item $p=1$ elements (uniform)
    \item Combination of $p=1$ on faces and edges with $p=3$ on volume (mixed).
    \item $p=2$ elements (uniform).
\end{enumerate}
We want to address that despite the de-facto agreement established in the community that $p=2$ offers the best compromise (mainly on synthetic models), other non-uniform schemes at different entity levels might be more advantageous for both synthetic and experimental models. Nevertheless, we note that each refinement combination performance depends on mesh quality (adapted $h-$refinement) and the input model. However, our numerical experiments provide key information to guide the end-user in analyzing the advantages and disadvantages of using a particular refinement scheme. 

Finally, we believe that our open-source 3D EM modeler and numerical experiments regarding supervised $h+p$ refinement using HEFEM will prove useful for the EM community. Hopefully, this analysis will provide an in-depth analysis of EM modeling's pros and cons when HEFEM and $h+p$ refinement are employed.

Our future research aims to implement the numerical schemes presented here in a fully automatic fashion (e.g., based on artificial intelligence, or introducing classical strategies of solve+estimate+refine). Also, we intend to perform simulations for unconstrained 3D anisotropic CSEM setups.

\section*{Code availability}
\label{code_availability}
The \petgem code is freely available at the home page (\href{petgem.bsc.es/}{\textit{petgem.bsc.es}}), at the PyPI repository (\href{pypi.org/project/petgem/}{\textit{pypi.org/project/petgem}}), at the GitHub site (\href{github.com/ocastilloreyes/petgem}{\textit{github.com/ocastilloreyes/petgem}}), or by requesting the author (\href{octavio.castillo@bsc.es}{\textit{octavio.castillo@bsc.es}}, \href{ocastilloreyes@gmail.com}{\textit{ocastilloreyes@gmail.com}}). In all cases, the code is supplied to ease the immediate execution on Linux platforms. User's manual and technical documentation (developer's guide) are provided in the \petgem archive.

\begin{acknowledgments}
The work of O.C-R. has received funding from "HPC Technology Innovation Lab, a Barcelona Supercomputing Center and Huawei research cooperation agreement (2020). O.C-R. has been also 65\% cofinanced by the European Regional Development Fund (ERDF) through the Interreg V-A Spain-France-Andorra program (POCTEFA2014-2020). POCTEFA aims to reinforce the economic and social integration of the French-Spanish-Andorran border. Its support is focused on developing economic, social and environmental cross-border activities through joint strategies favouring sustainable territorial development. Also, O.C-R. has received funding from the \emph{European Union's Horizon 2020 programme} under the \emph{Marie Sklodowska-Curie} grant agreement N$^\circ$ 777778. Furthermore, the development of \petgem has received funding from the \emph{European Union's Horizon 2020 programme}, grant agreement N$^\circ$~828947, and from the Mexican Department of Energy, CONACYT-SENER Hidrocarburos grant agreement N$^\circ$ B-S-69926. Finally, the authors would also like to thank the support of the \textit{Ministerio de Educaci\'on y Ciencia} (Spain) under Projects TEC2016-80386-P.

\end{acknowledgments}

\bibliographystyle{gji}
\bibliography{references}

\label{lastpage}

\end{document}